\input epsfx.tex
\input amstex
\catcode`@=12

\frenchspacing
\documentstyle{amsppt}
\magnification=\magstep1
\baselineskip=14pt
\vsize=18.5cm
\footline{\hfill\sevenrm version 20100506}
\def\C{\bold C}
\def\H{\bold H}
\def\R{\bold R}
\def\Q{\bold Q}
\def\Z{\bold Z}

\def\O{{\Cal O}}
\def\ga{{\goth a}}
\def\gb{{\goth b}}
\def\Pic{\text{\rm Pic}}
\def\Mat{\text{\rm Mat}}

\def\PSp{\text{\rm PSp}}
\def\Sp{\text{\rm Sp}}
\def\GL{\text{\rm GL}}
\def\PSL{\text{\rm PSL}}
\def\SL{\text{\rm SL}}
\def\Im{\text{\rm Im}}
\def\Tr{\text{\rm Tr}}

\def\G2#1{{\Gamma^{(2)}#1}}
\def\Gs2#1{{\Gamma_0^{(2)}#1}}
\def\Ys2#1{{Y_0^{(2)}#1}}
\gdef\Kronecker#1#2{\left( {#1 \over #2} \right)}
\def\endproof{\hfill$\square$\par}

\gdef\mapright#1{\ \smash{\mathop{\longrightarrow}\limits^{#1}}\ }
\gdef\hooklongrightarrow{\lhook\joinrel\longrightarrow}
\gdef\maptwohead#1{\ \smash{\mathop{\relbar\joinrel\twoheadrightarrow}\limits^{#1}}\ }

\newcount\refCount
\def\newref#1 {\advance\refCount by 1
\expandafter\edef\csname#1\endcsname{\the\refCount}}
\newref Bai   
\newref Co    
\newref CoA   
\newref CoB   
\newref Con   
\newref Del   
\newref EZ    
\newref Got   
\newref HW    
\newref IgB   
\newref IgA   
\newref IgC   
\newref Kob   
\newref Koe   
\newref Koh   
\newref KZ    
\newref Len   
\newref LenPom 
\newref Mah   
\newref RS    
\newref Ser   
\newref Shi   
\newref vW    
\newref Wash  
\newref Weng  

\topmatter
\title
Evaluating Igusa functions
\endtitle
\author Reinier Br\"oker, Kristin Lauter
\endauthor
\address
Brown University, Box 1917, 151 Thayer Street, Providence, RI 02912, USA
\endaddress
\address
Microsoft Research, One Microsoft Way, Redmond, WA 98052, USA
\email reinier\@math.brown.edu, klauter\@microsoft.com \endemail
\endaddress
\abstract
\hfuzz=3pt
The moduli space of principally polarized abelian surfaces
is parametrized by three Igusa functions. In this article we investigate a new
way to evaluate these functions by using Siegel Eisenstein series. We explain
how to compute the Fourier coefficients of certain Siegel modular forms using
classical modular forms of half-integral weight. One of the results in this
paper is an explicit algorithm to evaluate the Igusa functions to a
prescribed precision.
\endabstract
\endtopmatter

\document

\head 1. Introduction
\endhead
\noindent
The classical theory of complex multiplication gives an explicit description of
the Hilbert class field of an imaginary quadratic field: for a fundamental
discriminant $D<0$, the Hilbert class field of $K = \Q(\sqrt{D})$ is obtained
by adjoining the value $j((D + \sqrt{D})/2)$ to~$K$. Here, $j: \H \rightarrow
\C$ is the classical modular function with Fourier expansion $j(z) = 1/q+744+
196884q+\ldots$ in $q = \exp(2 \pi i z)$. There are various ways to compute the
minimal polynomial of $j((D + \sqrt{D})/2)$, and one of the most frequently
used approaches proceeds by evaluating the $j$-function to high precision.

The $j$-function is invariant under the action of $\SL_2(\Z)$ on
the upper half plane~$\H$. To evaluate $j(\tau)$, we may assume that 
$\tau$ is in the `standard' fundamental domain for $\SL_2(\Z)\backslash \H$
as described in e.g.\ [\Ser, Sec.\ VII.1.1]. 
The naive approach to evaluate $j(\tau)$ is to simply compute enough Fourier 
coefficients using for instance the
recursive formulas given in~[\Mah]. Alternatively, one can use the relation
$$
j(z) = 1728 {g_2(z)^3 \over g_2(z)^3 - 27g_3(z)^2} \eqno (1.1)
$$
expressing the $j$-function in terms of the {\it normalized Eisenstein\/}
series $g_2, g_3$ of weight $4$ and $6$.
Better results can be obtained [\Bai] by
using the Dedekind $\eta$-function defined by
$
\eta(z) = q^{1/24} \prod_{n=1}^{\infty}(1-q^n),
$
and which satisfies
$$
j(z) = \left( {(\eta(z/2)/\eta(z))^{24} + 16 \over (\eta(z/2)/\eta(z))^8}
       \right)^3.
$$
The sparsity of the $q$-expansion of the $\eta$-function makes it very
efficient for explicit computations.

The $j$-function is intrinsically linked to the theory of elliptic curves,
and the situation outlined above can be viewed as the `1-dimensional' case of
complex multiplication theory. In dimension 2, suitably chosen invariants of
principally polarized abelian surfaces generate abelian extensions of degree
4 CM-fields, see [\Shi, Sec.\ 15] for a precise statement. A popular choice of
invariants are the three {\it Igusa functions\/} $j_1,j_2,j_3$ defined
below. Just as evaluating the elliptic $j$-function has applications to 
elliptic curve cryptography, evaluating Igusa functions is an important
step in construction genus 2 curves suitable for use in cryptography, see
e.g.\ [\Weng].

The explicit evaluation of Igusa functions is less developed than its
dimension-1 counterpart. Most people use $\theta$-functions to evaluate 
Igusa functions. The (rather unwieldy) formulas expressing 
Igusa functions in terms of $\theta$-functions are given in 
e.g.\ [\Weng, pp.\ 441--442]. There is also a direct analogue of
formula~(1.1) which expresses the Igusa functions as rational functions in
the {\it Siegel Eisenstein series\/}~$E_w$. Indeed, Igusa [\IgA, p.~195] defines the normalized cusp forms
$$
\chi_{10} = -\frac{43867}{ 2^{12} \cdot 3^{5} \cdot 5^{2} \cdot 7 \cdot 53} (E_4 E_6 -E_{10})$$ and $$\chi_{12} = \frac{131\cdot 593}{2^{13} \cdot 3^{7} \cdot 5^{3} \cdot 7^{2} \cdot 337}
(3^2\cdot 7^2 E_4^3 +2 \cdot 5^3 E_4^6 -691 E_{12}).$$
With that, we have the three {\it Igusa functions}
$$ j_1 = 2\cdot 3^5 {\chi_{12}^5 \over \chi_{10}^6}, \quad j_2 = 2^{-3} 3^3
{E_4 \chi_{12}^3 \over \chi_{10}^4}, \quad j_3 = 2^{-5} \cdot 3 {E_6
\chi_{12}^2 \over \chi_{10}^3} + 2^{-3}\cdot 3^2 {E_4 \chi_{12}^3 \over
\chi_{10}^4} \eqno(1.2).$$
Igusa shows the equivalence with the definition of these functions 
in terms of theta functions in [\IgB, p. 848].
The analogue of the denominator $\Delta=g_2^3-27g_3^2$ appearing in~(1.1) is
the form~$\chi_{10}$. The form $\Delta$ is a classical cusp form of weight~12 and
$\chi_{10}$ is a Siegel cusp form of weight~10.

A mathematically natural question is whether we can use formula~(1.2) 
directly to evaluate
the Igusa functions, thereby bypassing the $\theta$-functions. The main
focus of this paper is to give an explicit algorithm to evaluate the Siegel
modular forms occuring in~(1.2) to some prescribed accuracy.
Our result gives a relatively easy way to analyze the precision necessary for 
the computation to succeed, and we give a {\it rigorous\/} complexity analysis 
for our method, something which has not been done for other approaches.

Although the asymptotic convergence of our algorithm is slower than the 
algorithm using theta functions, our approach has the advantage that there are 
fewer high precision multiplications required in the evaluation, and thus less 
precision loss and fewer rounding errors occur. 
Furthermore, we give a detailed analysis of the Eisenstein series and 
cusp forms, including an algorithm for computing them using classical 
modular forms of half-integral weight and explicit bounds on the size of the 
coefficients in their Fourier expansions. Indeed, one of the main 
contributions of the paper is the detailed analysis of various aspects of the 
computation of Siegel modular forms. Finally, our approach may lend itself to 
improvement in various ways and is a new direction in this area which could 
produce further progress.

Any Siegel modular form $f$ admits a Fourier expansion
$$
f(\tau) = \sum_T a(T) \exp(2\pi i\, \Tr(T\tau)) \eqno (1.3)
$$
where $T$ ranges over certain $2\times 2$-matrices with coefficients in
${1 \over 2}\Z$.
We propose to evaluate the functions occuring in~(1.2) by truncating the
sum in (1.3) to
only include matrices with trace below some bound.
The Eisenstein series are Siegel modular forms with
a considerable amount of extra structure. We show that computing the Fourier
coefficients of the Eisenstein series ultimately boils down to computing
Fourier coefficients of classical modular forms of half-integral weight.
One of the main results of this paper is the following theorem, proved in
Section~4.
\proclaim{Theorem 1.1}{For $A,C \in \Z_{\geq 0}$ and $B \in \Z$ with $B^2
\leq 4AC$, the Fourier coefficients of the Siegel Eisenstein
series $E_w$ for all matrices $\Bigl( {a \atop b/2}\thinspace {b/2 \atop c}
\Bigr)$ satisfying $0\leq a \leq A$, $0 \leq c \leq C$, $|b| \leq B$ can
be computed in time $O((ABC)^{1+\varepsilon})$ for every $\varepsilon>0$. The
constant in the $O$-symbol depends on the weight~$w$.}
\endproclaim\par
\ \par\noindent
By examining the {\it size\/} of the Fourier coefficients more closely, we
derive the following result in Section~6.
\proclaim{Theorem 1.2}{Let $\tau\in\H_2$ be given, and let $\delta =
\delta(\tau)$
be the supremum of all $\delta'\in\R$ such that $\Im(\tau) - \delta'1_2$ is
positive semi-definite. Assume that $\delta(\tau)\geq 1$. Assume
$\chi_{10}(\tau)$ is non-zero and choose $n \in \Z$ 
such that $|\chi_{10}(\tau)| \geq 10^{-n}$ holds.

For a positive integer~$k$, let $B \in \Z_{>0}$ be
such that
$$
\int_{B-1}^\infty 524093 t^{15} \exp(-2\pi t \delta(\tau)) \hbox{d}t\leq 10^{-k-\max\{22,6n\}}\eqno(1.4)
$$
holds.

Then the following holds: if we approximate the modular forms
$E_4, E_6, \chi_{10}, \chi_{12}$
using their truncated Fourier expansions consisting of all the matrices of
trace at most~$B$, then the values $j_1(\tau)$, $j_2(\tau)$,
$j_3(\tau)$ computed via the formulas in (1.2) are accurate to
precision~$10^{-k}$.}
\endproclaim
The condition $\delta(\tau) \geq 1$ is mostly for esthetic reasons. The proof
of Theorem~1.2, given in Section~6, readily gives a method to find $B$ in case
$\delta(\tau) < 1$. We assume in Theorem~1.2 that we can bound
$|\chi_{10}(\tau)|$ from below. This lower bound will allow us to bound the
precision loss that occurs when we divide by $\chi_{10}(\tau)$. Using the
explicit bounds on the Fourier coefficients of
$\chi_{10}$, proved in Section 5, we give a simple method to find a value of
$n$ in Section~6. This method works in general and does not depend on the
value of $\delta(\tau)$. Hence, Theorem~1.2 gives an effective method to
evaluate the three Igusa functions up to some prescribed precision.

Just as the elliptic $j$-function is invariant under $\SL_2(\Z)$, the Igusa
functions $j_1,j_2,j_3$ are invariant under the symplectic group~$\Sp_4(\Z)$.
Hence, we may translate the argument $\tau$ by a matrix $M \in \Sp_4(\Z)$ to
obtain an $\Sp_4(\Z)$-equivalent $\tau'\in \H_2$. The value $\delta(\tau')$
can be significantly different from $\delta(\tau)$, see e.g.\ Example 7.1. 
Before applying Theorem~1.2, we therefore move, using e.g.\ the method
from~[\vW], $\tau$ to the `standard' fundamental domain for 
$\Sp_4(\Z)\backslash\H_2$ described in~[\Got].

The outline of the article is as follows. In Section 2 we recall basic facts
about Siegel modular forms and their Fourier expansions. Section 3 introduces
Jacobi forms and their relation to Eisenstein series. The approach we follow
in this section is `classical' and most likely well-known to experts working
with Siegel modular forms. In Section 4 we go one step further, and relate
Jacobi forms to classical modular forms of half-integral weight. This gives
a very efficient method of computing the Fourier coefficients of
the 2-dimensional Eisenstein series. The functions $\chi_{10}$ and $\chi_{12}$
are Siegel cusp forms, and we explain in Section 5 how to compute the
Fourier coefficients of these forms.  We investigate the
convergence of the Fourier expansions of $E_4,E_6,\chi_{10}$ and
$\chi_{12}$ in Section~6.
This leads to the proof of Theorem~1.2. A final Section~7 contains two
detailed examples.

\head 2. Siegel modular forms
\endhead
\noindent
Let $\H_2 = \{ \tau \in \Mat_2(\C) \mid \tau = \tau^T, \Im(\tau) > 0 \}$ be
the Siegel upper half plane. With $J = \Bigl( {0 \atop -1_2} \thinspace
{1_2 \atop 0} \Bigr)$, the symplectic group $\Sp_4(\R)$ is defined as
$\Sp_4(\R) = \{ M \in \GL_4(\R) | M J M^T = J\}$. The group $\Sp_4(\R)$
naturally acts on the Siegel upper half plane via
$$
\Bigl( {a \atop c} \thinspace {b \atop d} \Bigr)\tau =
{a\tau + b \over c\tau+d},
$$
where dividing by $c\tau +d$ means multiplying on the right with the
multiplicative inverse of the $2\times 2$-matrix $c\tau +d$. The matrix
$-1_2$ acts trivially, and it is well
known that the automorphism group of $\H_2$ equals $\PSp_4(\R) =
\Sp_4(\R) / \{ \pm 1_2 \}$.

A holomorphic function $f: \H_2 \rightarrow \C$ is called a {\it Siegel
modular form\/} of weight $w \geq 0$ if it satisfies
$$
f( \Bigl( {a \atop c} \thinspace {b \atop d} \Bigr)\tau ) = \det(c\tau + d)^w
 f(\tau)
$$
for all $\tau$ and all matrices in the subgroup $\Sp_4(\Z) \subset \Sp_4(\R)$.
The integer $w$ is called the {\it weight\/} of the form~$f$. Whereas we
have to demand that $f$ is `holomorphic at infinity' for classical modular
forms $\H \rightarrow \C$, this is not necessary for Siegel modular forms.
Indeed, the {\it Koecher principle\/} implies that $f$ is bounded on sets
of the form $\{ \tau \in \H_2 \mid \Im(\tau) > \alpha 1_2 \}$ for $\alpha>0$,
see~[\Koe].

The matrix $\Bigl( {1 \atop 0} \thinspace {1 \atop 1} \Bigr)$ is contained
in $\Sp_4(\Z)$, and a Siegel modular function $f$ is invariant under the
transformation $\tau \mapsto \tau + 1$. In particular, a Siegel modular
function $f$ admits a {\it Fourier expansion\/}
$$
f(\tau) = \sum_T a(T) \exp(2\pi i \Tr(T\tau)).
$$
Here, the sum ranges over all symmetric matrices $T \in \Mat({1 \over 2}\Z)$
with integer diagonal entries. The coefficients $a(T)$ are called the
{\it Fourier coefficients\/} of $f$. By the Koecher principle, they are
zero in case $T$ is negative definite.

We embed the group $\GL_2(\Z)$ in $\Sp_4(\Z)$ via $M \mapsto \Bigl( {M \atop 0}
\thinspace {0 \atop (M^T)^{-1}} \Bigr)$. As $M^T$ has determinant $\pm 1$,
we see that a Siegel modular function $f$ is invariant under the transformation
$\tau \mapsto M \tau M^T$ for $M \in \GL_2(\Z)$. This invariance is the
key ingredient in the proof of the following well known lemma.

\noindent
\proclaim{Lemma 2.1}{The Fourier coefficients $a(T)$ of a Siegel modular form
$f$ satisfy $a(M^T T M) = a(T)$ for every $M \in \GL_2(\Z)$.}
\endproclaim
\noindent
{\bf Proof.} Writing $\tau = x + iy$ with $x,y \in \Mat_2(\R)$, the
Fourier coefficient $a(T)$ is given by
$$
a(T) = \int f(\tau) e^{-2 \pi i \Tr(T\tau)} \text{\rm d}x.
$$
Here, $\text{\rm d}x$ means the Euclidean volume of the space of $x$-coordinates
and the integral ranges over the `box' $-1/2 \leq x_{ij} \leq 1/2$. Using
the invariance of $f$ we compute
$$
a(M^T T M) = \int f(M \tau M^T) e^{-2 \pi i \Tr(T \thinspace M \tau M^T)}
\text{\rm d}x,
$$
and the lemma follows. \endproof\par
\ \par\noindent
In the
remainder of this section we investigate how many different values $a(T)$
attains for a fixed value of $n = \det(T)>0$ and a fixed Siegel modular
form~$f$.

To a matrix $T = \Bigl( {a \atop b/2} \thinspace {b/2 \atop c} \Bigr)$ with
$a,b,c \in \Z$ we associate the binary quadratic form
$f_T = aX^2 + bXY + cY^2$ of discriminant $b^2 - 4ac = -4n$. An explicit
check shows that for $M = \Bigl( {\alpha\atop \gamma} \thinspace {\beta\atop
\delta} \Bigr)$ the quadratic forms associated to $M^T T M$ equals
$$
f_{M^T T M} = f_T(\alpha X + \beta Y, \gamma X + \delta Y),
$$
which means that the $\GL_2(\Z)$-action on $\H_2$ is compatible with the
$\GL_2(\Z)$-action on quadratic forms of discriminant $-4n$. In fact, the
$\GL_2(\Z)$-action on quadratic forms originally considered by Lagrange is not
used that much as it leads to a `wrong' kind of equivalence. For quadratic
forms, the `correct' action is the action of the subgroup $\SL_2(\Z) \subset
\GL_2(\Z)$ studied by Gau\ss. The difference between these two actions
is implicit in the following lemma.

\noindent
\proclaim{Lemma 2.2}{Fix a Siegel modular form~$f$ and $n \in
{1\over 4}\Z_{> 0}$. Suppose that $-4n$ is a
fundamental discriminant and
let $\O$ be the maximal order of $\Q(\sqrt{-n})$. Then the set
$\{ a(T) \mid \det(T) = n \}$ has size at most
${1 \over 2} (\#\Pic(\O) + \#\{\ga\in\Pic(\O) \mid 2\ga = 0\})$.}
\endproclaim
\noindent
{\bf Proof.} If $-4n$ is fundamental, then any integer binary quadratic form
$aX^2 + bXY + cY^2$ of discriminant $-4n$ is {\it primitive\/}. The set of
$\PSL_2(\Z)$-equivalence classes of primitive quadratic forms of discriminant
$-4n$ is in bijection with the class group $\Pic(\O)$ via
$aX^2+bXY+cY^2 \mapsto a\Z + {-b+\sqrt{-4n}\over 2}\Z$ by [\Co, Th.\ 5.2.8].

It remains to investige when a $\GL_2(\Z)$-equivalence class decomposes as
2 disjunct $\SL_2(\Z)$-equivalence classes. If a fractional $\O$-ideal
$\ga$ is $\GL_2(\Z)$-equivalent but not $\SL_2(\Z)$-equivalent to $\gb$, then
$\gb$ equals the inverse $\ga^{-1}$ and we have $2\ga\not = 0$. The
lemma follows. \endproof\par
\ \par\noindent
For the general case of not necessarily fundamental discriminants, we note
that any binary quadratic form $aX^2+bXY+cY^2$ of discriminant $-4n$
determines a {\it primitive\/} quadratic form $(aX^2 + b XY + cY^2) /
\gcd(a,b,c)$ of discriminant $-4n/\gcd(a,b,c)^2$. Arguing as in the proof
of Lemma 2.2, we see that the set $\{ a(T) | \det(T) = n \in
{1\over 4}\Z_{>0}\}$ has at most
$$
{1 \over 2} \sum_{\O} \#\Pic(\O) + \#\{\ga\in \Pic(\O) \mid 2\ga = 0\}
$$
elements. Here, the sum ranges over all imaginary quadratic orders $\O$ that
contain the order of discriminant $-4n$.

\noindent
\proclaim{Corollary 2.3}{Let $m$ be the index of the order of discriminant
$-4n$ in the maximal order of the quadratic field $\Q(\sqrt{-n})$ and
let $\varphi$ denote
the Euler $\varphi$-function. For a fixed Siegel modular form~$f$, the set
$\{ a(T) | \det(T) = n \in {1 \over 4}\Z_{\geq0}\}$ then has as most
$2 \sqrt{n} \log(4n) (m/\varphi(m))^2$ elements. }
\endproclaim
\noindent
{\bf Proof.} The class number for the imaginary order of discriminant $D$
is bounded by $|D|^{1/2} \log |D|$ by [\LenPom, Sec.\ 2]. The result now follows
from the class number formula, see e.g.\ [\Len, Sec.\ 1.6].
\endproof

\head 3. Eisenstein series
\endhead
\noindent
For $w \geq 0$, the space $M_w$ of Siegel modular forms of weight~$w$ has a
natural structure of a
$\C$-vector space. For even $w \geq 4$, the primordial example of a degree~$w$
Siegel modular form is the Eisenstein series $E_w$ defined by
$$
E_w(\tau) = \sum_{c,d} (c\tau + d)^{-w}. \eqno(3.1)
$$
Here, the sum ranges over all inequivalent bottom
rows $(c \quad d)$ of elements of $\Sp_4(\Z)$ with respect to
left-multiplication by~$\SL(2,\Z)$.
The restriction $w \geq 4$ comes from the fact that the
expression in (3.1) does not converge for $w=2$.

The direct product $M = \coprod_{w=0}^{\infty} M_w$ has a natural structure
of a graded $\C$-algebra. By restricting the product to {\it even\/} $w$,
we get a graded subalgebra~$M^e$. The following lemma gives the structure
of these two algebras.

\proclaim{Lemma 3.1}{The Eisenstein series $E_4$, $E_6$, $E_{10}$ and
$E_{12}$ are algebraically independent and generate~$M^e$. There exists
a polynomial $P$ in $4$ variables such that $M$ is isomorphic
to $M^e[X]/(X^2-P(E_4,E_6,E_{10},E_{12}))$. The element $\overline X$
corresponds to a Siegel modular form of weight~$35$.}
\endproclaim
\noindent
{\bf Proof.} The first statement can be found in [\IgA, pp.\ 194--195]. The
second statement is proven in [\IgB] with an explicit polynomial~$P$ at
page~849. \endproof
\ \par\noindent
The remainder of this section is devoted to deriving a `formula' for the
Fourier coefficient $a(T)$ of the Eisenstein series $E_w$. The approach we
follow is intrinsically related to the theory of {\it Jacobi forms\/}, see [\EZ]
for a good introduction. Let $f : \H_2 \rightarrow \C$ be a Siegel modular
form of weight~$w$. We write $\tau \in \H_2$ as $\tau = \Bigl( {\tau_1 \atop
\varepsilon} \thinspace {\varepsilon \atop \tau_2}\Bigr)$. Because $f$ is
periodic with respect to $\tau_2$, it admits a Fourier expansion
$$
f(\tau) = \sum_{m=0}^{\infty} \varphi_m(\tau_1,\varepsilon) e^{2 \pi i m\tau_2}
$$
where $\varphi_m$ is a function from $\H \times \C$ to $\C$. The functions
$\varphi_m$ have the following properties:\par
\medskip
\item{$\diamond$}{\quad$\varphi_m({a\tau_1+b \over c\tau_1+d},
{\varepsilon \over c \tau_1 + d}) = (c \tau_1 + d)^{\omega} e^{2\pi imc
\varepsilon/(c \tau_1+d)} \varphi_m(\tau_1,\varepsilon), \qquad
\Bigl( {a\atop c} \thinspace {b \atop d} \Bigr) \in \SL_2(\Z)$}
\smallskip
\item{$\diamond$}{\quad$\varphi_m(\tau,\varepsilon+\lambda\tau+\mu) =
e^{-2\pi im(\lambda^2\tau_1+2 \lambda\varepsilon)}
\varphi_m(\tau_1,\varepsilon), \qquad (\lambda,\mu)\in\Z^2$}
\smallskip
\item{$\diamond$}{\quad $\displaystyle\varphi_m \hbox{\ admits a Fourier expansion of the
form\ } \sum_{n=0}^{\infty}\sum_{r \in \Z \atop {\scriptscriptstyle r^2 \leq
4nm}} c(n,r) e^{2\pi i(n\tau_1+r\varepsilon)}.$}
\medskip\noindent
The first two properties follow from the transformation law of Siegel modular
forms under the symplectic matrices\par
\medskip
\centerline{
$
\pmatrix
a&0&b&0\\ 0&1&0&0\\ c&0&d&0\\ 0&0&0&1\\
\endpmatrix
$
\qquad and \qquad
$
\pmatrix
1&0&0&\mu\\ \lambda&1&\mu&0\\ 0&0&1&-\lambda\\ 0&0&0&1\\
\endpmatrix
$
}\par
\medskip\noindent
and the third property follows from the Koecher principle.

A holomorphic function $g: \H \times \C \rightarrow \C$ satisfying the three
properties above for some $w$ and $m$ is called a {\it Jacobi form\/} of
weight~$w$ and index~$m$.
Jacobi forms can be seen as an `intermediate' between Siegel modular forms
and classical modular forms. Indeed, the `Fourier coefficients' of a Siegel
modular form of weight~$w$ are Jacobi forms of weight~$w$ and for a Jacobi
form $g$, the function $g(\tau,0)$ is a classical modular form of weight~$w$.

The space of all Jacobi forms of weight $w$ and index $m$ is denoted by
$J_{w,m}$, and we have maps
$$
M_w \hooklongrightarrow \prod_{m\geq 0} J_{w,m} \maptwohead{\text{\rm pr}} J_{w,1},
$$
where pr denotes the projection onto the first factor.
For this article, the key property of Jacobi forms is that we can also
construct a map $J_{w,1} \rightarrow M_w$ which will allow us to identify
certain Siegel modular forms with its `first' Jacobi form. As we have
$J_{w,1} = 0$ for {\it odd\/}~$w$ by~[\EZ, Th.~2.2], we restrict to
{\it even\/} weight~$w$ for the remainder of this section.

For $m\geq 0$,
we define the `Hecke operator' $V_m: J_{w,1} \rightarrow J_{w,m}$ as follows.
For $g \in J_{w,1}$ with Fourier expansion $\sum_{n,r} c(n,r)
e^{2 \pi i(n\tau_1 +r\varepsilon)}$, we put
$$
V_m(g) = \sum_{n,r} \left( \sum_{a \mid \gcd(n,r,m)} a^{w-1} c\left(
{nm\over a^2}, {r\over a}\right) \right) e^{2\pi i(n\tau_1+r\varepsilon)}
$$
for $m>0$. This is the natural generalization of the Hecke operators for
classical modular forms, see e.g.\ [\Ser, Prop.\ VII.12]. For $m=0$, we put
$$
V_0(g) = -{B_w c(0,0)\over 2w}\left(1-{2w \over B_w}\sum_{n\geq 1} \sigma_{w-1}(n)
e^{2\pi i n\tau_1}\right)
$$
with $\sigma_n(x)$ the sum of the $n$th powers of the divisors of~$x$ and
$B_w$ the $w$th Bernoulli number defined by $t/(e^t - 1)=\sum_{n=0}^\infty
B_n t^n/n!$. In
particular, the function $V_0(g)$ is a multiple of the classical Eisenstein
series of weight~$w$.
It is not hard to show that the function
$$
\Psi(g) = \sum_{m \geq 0} V_m(g)(\tau_1,\varepsilon) e^{2\pi i m\tau_2}
$$
defines a Siegel modular form of weight $w$, see [\EZ, Th.\ 6.2].

\proclaim{Lemma 3.2}{The map $\Psi: J_{w,1} \rightarrow M_w$ is injective.}
\endproclaim
\noindent
{\bf Proof.} This follows directly from the fact that the composition
$$
J_{w,1} \mapright{\Psi} M_w \hooklongrightarrow \prod_{m\geq 0} J_{w,m}
\maptwohead{\text{\rm pr}} J_{w,1}
$$
is the identity. \endproof\par
\ \par\noindent
We stress that the map $\Psi$ is in general {\it not\/} surjective. The
image $\Psi(J_{k,1})$ is known as the {\it Maa\ss\ Spezialschar\/}. However,
the Eisenstein series $E_w \in M_w$ do occur at the image of a Jacobi form.
They are the images of the {\it Jacobi Eisenstein series\/} $E_{w}^J$
defined by the (rather awkward looking) formula $E_{w}^J(\tau,z)=$
$$
{1 \over 2} \sum_{c,d\in \Z \atop {\scriptscriptstyle
\gcd(c,d)=1}}\sum_{\lambda\in\Z} (c\tau+d)^{-w} \exp\left(2\pi i\left(
\lambda^2 {a\tau+b \over c\tau+d} + 2 \lambda {z \over c\tau+d}-{c z^2\over
c\tau+d} \right)\right)
$$
for $w \geq 4$. Here, $a$ and $b$ are integers such that $\Bigl( {a \atop c}
\thinspace {b \atop d} \Bigr)$ is contained in $\SL_2(\Z)$.

\proclaim{Lemma 3.3}{We have $E_w =\Psi\Bigl( {-2w \over B_w} E_{w}^J\Bigr)$.}
\endproclaim
\noindent
{\bf Proof.} It follows from [\EZ, Th.\ 6.3] that $E_w$ is a multiple of
$\Psi(E_{w}^J)$. Both the Siegel Eisenstein series $E_w$ and the Jacobi
Eisenstein series $E_w^J$ are normalized with constant coefficient~1. The
lemma follows.\endproof
\ \par\noindent
It is now a straightforward matter to compute the Fourier coefficients of
the Siegel Eisenstein series. The result is the following theorem.

\proclaim{Theorem 3.4}{Let $E_w$ be the Siegel Eisenstein series of weight~$w$,
and let $T = \Bigl( {a \atop b/2} \thinspace {b/2 \atop c} \Bigr) \in
\Mat({1 \over 2}\Z)$ be a positive semi-definite matrix with integer entries
on the diagonal. Write $D = b^2-4ac\leq 0$ and let $D_0$ be the
discriminant of $\Q(\sqrt{D})$. Then the Fourier coefficient $a(T)$ equals $1$
for $a=b=c=0$ and
$$
{-2w \over B_{w}} \sum_{d | \gcd(a,b,c)} d^{w-1} \alpha(D/d^2)
$$
otherwise. Here, $B_{k}$ is the $k$th Bernoulli number and $\alpha$ is defined by
$\alpha(0) = 1$ and
$$
\alpha(D) = {1 \over \zeta(3-2w)} C(w-1,D) \qquad\qquad (D < 0)
$$
where $C$ is Cohen's function defined by
$$
C(s-1,D) = L_{D_0}(2-s) \sum_{d \mid f} \mu(d) \Kronecker{D_0}{d} d^{s-2} \sigma_{2s-3}(f/d), \qquad\qquad D = D_0 f^2.
$$
Here, $\zeta$ denotes the Dedekind
$\zeta$-function, $L_{D_0}$ is the quadratic Dirichlet $L$-series,
$\mu$ is the M\"obius function, $\Kronecker{\cdot}{\cdot}$ is the Kronecker
symbol and $\sigma_n(x)$ denotes the sum of the $n$th powers of the
divisors of~$x$.}
\endproclaim

\noindent
{\bf Proof.} By~[\EZ, Th.\ 2.1], the term $\alpha(D')$ equals the Fourier
coefficient $\alpha(n,r)$ of the Jacobi Eisenstein series $E_{w}^J$ with
$D' = r^2-4n$. By Lemma 3.3, we have to apply the Hecke operators $V_m$ to
these coefficients. The theorem follows.\endproof

\ \par\noindent
A formula for $a(T)$ is also given in Corollary 2 to [\EZ, Th.\ 6.3]. In this
formula, the Bernoulli numbers and the $\zeta$-function from Theorem
3.4 are missing.

We see that Theorem 3.4 gives a much better bound than Lemma 2.2 for the
cardinality of $\{ a(T) \mid \det(T) = n \in {1 \over 4} \Z\}$ for the
Eisenstein series. Indeed, for
fundamental discriminants $-4n$, we have only {\it one\/} Fourier
coefficient~$a(T)$. In general, the number of coefficients is bounded by the
number of square divisors of $-4n$ which in turn is bounded by
$O(n^{\varepsilon})$ for all $\varepsilon>0$. These bounds hold in general
for functions in the Spezialschar $\Psi(J_{w,1}) \subset M_w$. Indeed, the
Fourier coefficients $c(n,r)$ of a function $g \in J_{w,1}$ only depend
on the value $4n-r^2$, cf.\ [\EZ, Th.\ 2.2].

\proclaim{Corollary 3.5}{Let $n_k$ be the numerator of the $k$th
Bernoulli number $B_k$. Then the Fourier coefficient $a(T)$ of the
Siegel Eisenstein series $E_w$ for the
matrix $T = \Bigl( {a \atop b/2} \thinspace {b/2 \atop c} \Bigr)$ is
contained in the set $1/(n_wn_{2w-2})\Z\subset\Q$.}
\endproclaim
\noindent
{\bf Proof.} As we have $\zeta(3-2w) = -B_{2w-2}/(2w-2)$ all we have to do
is examine the denominator of the value $L_{D_0}(2-w)$ occuring in
Theorem~3.4. This is most easily done using $p$-adic $L$-series as in
[\CoA, Ch.\ 11]. The corollary follows from [\CoA, Cor.\ 11.4.3]
except in the following case: the discriminant of $\Q(\sqrt{b^2-4ac})$
equals $-p$ for an odd prime $p$ with $w-1 \equiv (p-1)/2 \bmod (p-1)$. If
this is the case, we a priori find that the denominator of the $L$-value
could be divisible by $p(w-1)$.
However, the prime $p$ then satisfies $(p-1) \mid (2w-2)$ and
by the Clausen-von Staudt theorem [\CoA, Cor.\ 9.5.15] the prime $p$ also
occurs in the denominator of $B_{2w-2}$. Finally, $w-1$ is a divisor of
the denominator of $\zeta(3-2w)$. \endproof

\proclaim{Corollary 3.6}{The Fourier coefficient $a(T)$  of the
Siegel Eisenstein series $E_w$ for the
matrix $T = \Bigl( {a \atop b/2} \thinspace {b/2 \atop c} \Bigr)$
satisfies $|a(T)| = O((4ac-b^2)^{w-3/2})$ if $b^2-4ac$ is non-zero.}
\endproclaim
\noindent
{\bf Proof.} Using the functional equation for Dirichlet $L$-series, see
e.g.\ [\CoA, Th.\ 10.2.6], we bound $L_{D_0}(2-w) = O(D_0^{w-3/2})$. The
inequalities
$$
{\sigma_n(x) \over x^n} = \sum_{d \mid x} {1 \over d^n} \leq
\sum_{d=1}^{\infty} {1 \over d^n} = \zeta(n) < \infty
$$
give $\sigma_n(x) = O(x^n)$ for $n>1$. It follows that the $c(D')$ in
Theorem~3.4 is of size $O({D'}^{w-3/2})$. As $\sum_{d \mid n} d^{w-1} /
d^{2w-3}$ is finite for $w\geq 4$ and $n\rightarrow \infty$, the corollary
follows. \endproof
\proclaim{Remark 3.7}{It is not hard to make the constant $c$ in the $O$-symbol
explicit. One can take
$$
c = \left|{4w (w-2)! \, \zeta(w-1)^2 \zeta(2w-3)\zeta(w-2) \over
\pi^{w-1}\zeta(3-2w)B_w}\right|.
$$
}
\endproclaim

\head 4. Computing special values of $L$-series
\endhead
\noindent
The hard part in computing Fourier coefficients of Siegel Eisenstein series
is computing the special values of $L$-series occuring in Theorem~3.4. If
the discriminant of the quadratic field $\Q(\sqrt{b^2-4ac})$,
corresponding to the matrix $\Bigl( {a \atop b/2} \thinspace {b/2 \atop c}
\Bigr)$, is small these computations can be efficiently done employing
generalized Bernoulli numbers as we now explain.

For $n\geq 1$, we let $\chi_n$ be the quadratic Dirichlet character modulo $n$
and define the {\it $\chi_n$-Bernoulli numbers} $B_k(\chi_n)$ by the expansion
$$
{\sum_{r=1}^n \chi_n(r) t e^{rt} \over e^{nt}-1} = \sum_{k\geq 0}
{B_k(\chi_n) \over k!} t^k \in \Q[t]. \eqno(4.1)
$$
The generalized Bernoulli numbers $B_k(\chi_n)$ equal the ordinary
Bernoulli numbers $B_k$ for $n=1$ and $k \geq 2$.

\proclaim{Lemma 4.1}{For $n\geq 1$ and $w\geq 2$, we have $L_n(2-w) =
-B_{w-1}(\chi_n)/(w-1)$.}
\endproclaim
\noindent
{\bf Proof.} See [\Wash, Th.\ 4.2]. \endproof \par
\ \par\noindent
The values $B_{w-1}(\chi_n)$ can easily be computed using the definition
(4.1) for small $w$ and $n$. For evaluating the Igusa functions, we are
only interested in the values $w=4,6,10,12$ and by computing $B_{11}(\chi_n)$
we get the other values $B_{9}(\chi_n), B_5(\chi_n)$ and $B_3(\chi_n)$ `for
free'.

To compute the Fourier coefficients of the Eisenstein series $E_w$ for
{\it large\/} values of~$D = b^2-4ac$, we clearly need another method. It is
a relatively well-known fact that Jacobi forms of even weight and index 1
`correspond to' classical modular forms of half-integral weight.
Explicitly, for a discriminant $D<0$, we define
\medskip
\halign{\qquad\qquad$#$ & $\displaystyle{#}$\cr
\alpha_w(D) & = {1 \over \zeta(3-2w)} C(w-1,D)\cr
&  = {1 \over \zeta(3-2w)} L_{D_0}(2-w)
\sum_{d \mid f} \mu(d) \Kronecker{D_0}{d} d^{w-2} \sigma_{2w-3}(f/d)\cr}
\medskip\noindent
as in Theorem 3.4. Here, $D_0$ is the discriminant of the quadratic field
$\Q(\sqrt{D})$ and $f$ satisfies $D_0f^2 = D$. We put $\alpha_w(0)=1$, and
$\alpha_w(D)=0$ if $D<0$ is not a discriminant.
The function $H_w : \H \rightarrow \C$ defined by
$$
H_w(z) = \sum_{n=0}^{\infty} \alpha_w(-n) q^n \qquad (q = \exp(2\pi i z))
$$
is known as {\it Cohen's function\/}.\par
\proclaim{Lemma 4.2}{Let $H_w$ be defined as above. Then $H_w$ is a modular
form of weight~$w-1/2$ for the congruence subgroup $\Gamma_0(4)$.}
\endproclaim
\noindent
{\bf Proof.} See [\CoB, Th.\ 3.1], or an alternate proof in [\Kob, Prop.\ IV.6].\endproof\par
\ \par\noindent
\proclaim{Remark}{}The bound $\alpha_w(n) = O(n^{w-3/2})$ from the proof of
Corollary~3.6 is in nice accordance with the general result that the
Fourier coefficients of a modular form of weight~$k$ are of size $O(n^{k-1})$.
\endproclaim
\noindent
As the $\C$-vector space of modular forms of fixed (half-integral) weight
is finite dimensional, we can easily compute coefficients of $H_w$ given a
basis for the vector space. It is not hard to show that the function
$$
\theta(z) = \sum_{n \in \Z} q^{n^2} = 1 + 2 \sum_{n=1}^{\infty} q^{n^2}
\qquad (q = e^{2 \pi i z})
$$
is a modular form of weight~$1/2$ for $\Gamma_0(4)$. The function
$$
\widetilde\theta(z) = \theta^4(z+1/2) = \Bigl(1 + 2 \sum_{n=1}^{\infty}
(-1)^n q^{n^2} \Bigr)^4 \qquad (q = e^{2\pi i z})
$$
is therefore a modular form of weight~2. Analogous to the proof of
[\Kob, Prop.\ IV.4], it follows that $\theta$ and $\widetilde\theta$
generate the $\C$-{\it algebra\/} of all modular forms. The main advantage
of choosing this basis is that $\theta$ is very lacunary.

\proclaim{Proposition 4.3}{The following equalities hold: }
\medskip
\halign{\quad $#$ \thinspace & $#$ & \thinspace $\displaystyle #$ \cr
H_4 &=& {\theta^7 + 7\theta^3 \widetilde\theta\over 8}\cr
\noalign{\smallskip}
H_6 &=& {- \theta^{11} + 22\theta^7\widetilde\theta + 11\theta^3\widetilde
       \theta^2 \over 32}\cr
\noalign{\smallskip}
H_{10} &=& {-43867\theta^{19}+725876\theta^{15}\widetilde\theta+12824886
          \theta^{11}\widetilde\theta^2+8845412\theta^7
          \widetilde\theta^3+107597\theta^3\widetilde\theta^4 \over
          22459904}\cr
\noalign{\smallskip}
H_{12} &=& {77683\theta^{23}+212405\theta^{19}\widetilde\theta+38627902
          \theta^{15}\widetilde\theta^2+100820362\theta^{11}\widetilde
          \theta^3 \over 159094784} + \cr
\noalign{\smallskip}
       && +{19313951\theta^7\widetilde\theta^4+42481\theta^3\widetilde
           \theta^5 \over 159094784}. \cr}
\endproclaim
\noindent
{\bf Proof.} Using Lemma 4.1, we compute the first few Fourier coefficients
of $H_w$ for $w=4,6,10,12$. With the obervation that $H_w$ equals an {\it
isobaric\/} polynomial in $\theta$ and $\widetilde\theta$, we have to
solve a system of $w/2$ equations in $w/2$ unknowns. The theorem follows.
\endproof\par
\ \par\noindent
This theorem allows us to compute the first $N$ coefficients of $H_w$ in
time $O(N^{1+o(1)})$ using fast multiplication techniques. This leads
to the theorem stated in the introduction. An important conclusion is
that it is much faster to compute $L$-values {\it simultaneously\/} than to
compute them individually.

\proclaim{Corollary 4.4}{For $A,C \in \Z_{\geq 0}$ and $B \in \Z$ with $B^2
\leq 4AC$, the Fourier coefficients of the Siegel Eisenstein
series $E_w$ for all matrices $\Bigl( {a \atop b/2}\thinspace {b/2 \atop c}
\Bigr)$ satisfying $0\leq a \leq A$, $0 \leq c \leq C$, $|b| \leq B$ can
be computed in time $O((ABC)^{1+\varepsilon})$ for every $\varepsilon>0$. The
constant in the $O$-symbol depends on the weight~$w$.}
\endproclaim

\head 5. Cusp forms
\endhead
\noindent
The techniques explained in Sections 3 and 4 allow us to efficiently compute
the Fourier coefficients of Siegel Eisenstein series. This suffices for
evaluating Igusa functions, since these functions are rational expressions
in $E_w$ for $w=4,6,10,12$. However, the denominators of the Igusa functions
have more structure: they are Siegel cusp forms. It is a natural question
to ask if we can compute the Fourier coefficients of $\chi_{10}$
{\it directly\/} via Jacobi forms. We explain this method in this Section.

Let $M^1_w$ be the vector space of classical modular forms of integral
weight~$w$,
and let $M^1 = \coprod_{w\geq 0} M^1_w$ be the space of all classical
modular forms. It is well known that we have $M^1 \cong \C[E_4^1,E_6^1]$,
with $E_w^1$ the classical Eisenstein series of weight~$w$, see
[\Ser, Cor.\ 2 to Th.\ VII.4]. We define the {\it Siegel operator\/}
$S: M \rightarrow M^1$ as follows. For a Siegel modular form
$f: \H_2 \rightarrow \C$ with Fourier expansion $f(\tau) = \sum_T a(T)
\exp(2\pi i \Tr(T\tau))$ we put
$$
S(f) = \sum_{n \geq 0} a\left( \Bigl( {n\atop 0} \thinspace
{0 \atop 0} \Bigr) \right) e^{2\pi i n\tau_1}, \qquad\qquad
\hbox{\ with\ } \tau = \Bigl( {\tau_1\atop \varepsilon} \thinspace
{\varepsilon\atop\tau_2} \Bigr).
$$
The Siegel operator is a ring homomorphism $M \rightarrow M^1$, and it
maps Eisenstein series to Eisenstein series. In fact, for the Eisenstein
series $E_w$, it is the composition of the maps
$$
M_w \mapright{} \prod_{m\geq 0} J_{w,m} \maptwohead{\text{\rm pr}}
J_{w,0} \mapright{} M_w^1,
$$
introduced in Section~2.

A Siegel modular form $f$ is called a {\it cusp form\/} if it satisfies
$S(f) = 0$. Equivalently, $f$ is a cusp form if and only if the Fourier
coefficients $a(T)$ are zero for all semi-definite $T$ that are not
definite. It follows from well-known identities between classical Eisenstein
series that
$$
\chi_{10} = -43867\cdot 2^{-12} \cdot 3^{-5} \cdot
5^{-2} \cdot 7^{-1} \cdot 53^{-1} (E_4 E_6 -E_{10})
$$
and
$$
\chi_{12} = 131\cdot 593\cdot 2^{-13} \cdot 3^{-7} \cdot 5^{-3} \cdot 7^{-2}
\cdot 337^{-1} (3^2\cdot 7^2 E_4^3 +2 \cdot 5^3 E_4^6 -691 E_{12}),
$$
are cusp forms. The constants
in $\chi_{10}$ and $\chi_{12}$ should be regarded as `normalization factors'.

\proclaim{Lemma 5.1}{The ideal of cusp forms in $M^e$ is generated by
$\chi_{10}$ and $\chi_{12}$. The ideal of cusp forms in $M$ is generated by
$\chi_{10}, \chi_{12}$ and a modular form $\chi_{35}$ of weight~35
corresponding to $\overline X$ in Lemma 3.1}
\endproclaim
\noindent
{\bf Proof.} See [\IgC, Th. 3].\endproof\par
\ \par\noindent
It is well-known that the cusp forms $\chi_{10}$ and $\chi_{12}$ are
contained in the Maa\ss\ Spezialschar~$\Psi(J_{k,1})$, the
gest of the proof being [\EZ, Th.\ 6.3]. A Jacobi form $g\in J_{w,m}$ is called
a {\it cusp form\/} if its Fourier coefficients $c(n,r)$ are zero for
$4nm-r^2=0$. In particular, the map
$$
M_w \rightarrow \prod_{m \geq 0} J_{w,m} \maptwohead{\text{\rm pr}} J_{w,1}
$$
maps Siegel cusp forms to Jacobi cusp forms.
In weight $10$ and $12$ we have the Jacobi cusp forms
$$
\varphi_{10,1} = {1 \over 144}(E_6^1 E_{4,1} - E_4^1 E_{6,1}) \qquad \hbox{and}
\qquad \varphi_{12,1} = {1 \over 144}((E_4^1)^2 E_{4,1}-E_6 E_{6,1}),
$$
with $E_4^1 = 1 + 240 \sum_{n>0} \sigma_3(n) q^n$ and $E_6^1 = 1 - 504
\sum_{n>0} \sigma_5(n) q^n$ the classical Eisenstein series. The factor 144
should again be regarded as a normalization factor.

\proclaim{Lemma 5.2}{We have $\chi_{10} = \Psi(-\varphi_{10,1}/4)$ and
$\chi_{12} = \Psi(\varphi_{12,1}/12)$.}
\endproclaim
\noindent
{\bf Proof.} The cusp forms $\chi_{10}$ and $\chi_{12}$ are contained in the
Spezialschar and therefore occur as images of Jacobi cusp forms. The
spaces of Jacobi cusp forms of weight $10$ and $12$ are 1-dimensional by
[\EZ, Th.\ 3.5]. Using Theorem 3.4, we compute
$$
a\left( \Bigl({1 \atop 0} \thinspace {0 \atop 1} \Bigr) \right) =
{1 \over 2}
$$
for a Fourier coefficient of $\chi_{10}$ and
$$
c(1,0) = -2
$$
for the corresponding coefficient of $\varphi_{10,1}$. The result for
$\chi_{10}$ follows. The computation for $\chi_{12}$ yields
$$
a \left( \Bigl({1 \atop 0} \thinspace {0 \atop 1} \Bigr) \right) =
{5 \over 6} \qquad \hbox{and} \qquad c(1,0) = 10.
$$
The lemma follows. \endproof\par
\ \par\noindent
To compute the Fourier coefficients of $\varphi_{10,1}$ and $\varphi_{12,1}$
we note that the coefficients $c_{\varphi_{10,1}}(n,r)$ and
$c_{\varphi_{12,1}}(n,r)$ only depend on the value of~$4n-r^2\geq 0$.
Furthermore, the functions
$$
K_{10} = \sum_{k\geq 0} c_{\varphi_{10,1}}(k) q^k \qquad \hbox{and} \qquad
K_{12} = \sum_{k\geq 0} c_{\varphi_{12,1}}(k) q^k
$$
are classical modular {\it cusp forms\/}, for the group $\Gamma_0(4)$, of
weight $9{1\over 2}$ and $11{1\over 2}$ respectively by [\EZ, Th.\ 5.4].
\proclaim{Proposition 5.3}{Let $\theta$ and $\widetilde\theta$ be as in Section~4.
Then we have
\medskip
\halign{\quad $#$ \thinspace & $#$ & \thinspace $\displaystyle #$ \cr
K_{10} &=& {\theta^{15}\widetilde\theta -3\theta^{11}\widetilde\theta^2 +
3\theta^7\widetilde\theta^3 - \theta^3\widetilde\theta^4\over 4096}\cr
\noalign{\smallskip}
K_{12} &=& {5\theta^{19}\widetilde\theta - 16\theta^{15}\widetilde\theta^2
+18\theta^{11}\widetilde\theta^3 - 8 \theta^7\widetilde\theta^4+\theta^3
\widetilde\theta^5 \over 16384}.\cr}
}
\endproclaim
\medskip\noindent
{\bf Proof.} Analagous to the proof of Theorem~4.3. \endproof
\par\ \par\noindent
It should come as no surprise that there are no terms $\theta^{19}$ and
$\theta^{23}$ occuring in Proposition~5.3. Indeed, the forms $K_{10}$ and
$K_{12}$ are cusp forms and therefore vanish at $q=0$ whereas $\theta$
does not vanish at $q=0$.

Proposition~5.3 allows us to evaluate the Siegel cusp forms $\chi_{10}$ and
$\chi_{12}$ at arbitrary points~$\tau\in\H_2$. For the proof of Theorem~1.2,
we need a bound on the size of the Fourier
coefficients of $\chi_{10}$ and $\chi_{12}$ as well. We need an {\it
explicit\/} bound, like the bound in Remark~3.7.

The `Resnikoff-Salda\~na conjecture' ([\RS])
$$
|a(T)| = O((\det T)^{w/2-3/4+\varepsilon})
$$
for the size of a Fourier coefficient $a(T)$ of a Siegel modular cusp form
of weight~$w$ is known to be false in general. At this moment, the best
known result is
$$
|a(T)| = O((\det T)^{w/2-13/36+\varepsilon})
$$ for every $\varepsilon>0$, see~[\Koh]. We will prove in the remainder of
this section that the Fourier coefficients of $g_{10}$ and $g_{12}$
satisfy
$$
|a(T)| = O((\det T)^{w/2-1/2+\varepsilon}),
$$
and we will make the constant in the $O$-symbol explicit. The reason that our
bound is better than the bound in~[\Koh] is that our Siegel modular forms
lie in the Maa\ss\ Spezialschar, and this allows us to give a stronger
bound. In fact, we will show that if the Lindel\"of-hypothesis is true, the
Fourier coefficients are of size~$O((\det T)^{w/2-3/4+\varepsilon})$.

First we will bound the Fourier coefficients of $K_{10}$ and $K_{12}$
explicitly. One approach would be to adapt `Hecke's proof' [\Ser, Th.\ VII.5]
for cusp forms. This technique would yield a bound of $O(n^{4.75})$ for
$g_{10}$, where we can make the constant in the $O$-symbol explicit.
However, our modular forms have considerably more structure and we will use a
variant of Waldspurger's formula to obtain a better bound.

The modular forms $K_{10}$ and $K_{12}$ have the property that their Fourier
coefficients $a(n)$ are zero for $n \equiv 1,2 \bmod 4$. As a consequence,
see~[\EZ, Sec\. 6], both functions are Hecke eigenforms. For every Hecke eigenform
$f$ of weight $w-1/2$, Shimura constructs, see e.g.\ [\EZ, Sec.\ 5], an integral weight cuspform $g$ of
weight $2w-2$ with the property that
$$
|a(D)|^2 = {\langle g,g \rangle \over \langle f,f \rangle} {(w-2)! \over \pi^{w-1}} L(g,\chi_D,w-1) |D|^{w-3/2}
\eqno(5.1)
$$
holds for the Fourier coefficient $a(D)$ of~$f$. In formula~(5.1), known as
Waldspurger's formula, $\langle \cdot,\cdot \rangle$ denotes the usual
Petersson inner product, and $L(g,\chi_D,s)$ is the $L$-series associated
to~$g$, twisted by the quadratic Dirichlet character~$\chi_D$.

\proclaim{Lemma 5.4}{Let $g_{18}, g_{22}$ be the modular forms associated to
$K_{10}, K_{12}$ respectively under Shimura's construction. Then we have
$$
{\langle g_{18},g_{18} \rangle \over \langle K_{10},K_{10}\rangle} \leq 75634 \qquad
\hbox{\ and\ }
{\langle g_{22},g_{22} \rangle \over \langle K_{12},K_{12}\rangle} \leq 1197339.
$$
}
\endproclaim
\noindent
{\bf Proof.} As the space of weight 18 cusp forms is one-dimensional, we have
$g_{18} = \Delta E_6^1$. Likewise, $g_{22} = \Delta E_{10}^1$. By formula~(5.1) we
have
$$
{\langle g_{18}, g_{18} \rangle \over \langle K_{10},K_{10} \rangle} =
{L(g_{18},\chi_D,9) \over |a(D)|^2} \cdot {8! \over \pi^9} |D|^{8.5}
$$
{\it for every\/} discriminant~$D$ for which the Fourier coefficient $a(D)$
of $K_{10}$ is nonzero. Since we can
compute the Fourier coefficients of $K_{10}$, it suffices to explicitly
evaluate the $L$-series at the center of the critical strip.

Since $g_{18}$ is a Hecke eigenform, the formula
$$
(2\pi)^{-s} \Gamma(s) L(g_{18},s) = \int_0^\infty g_{18}(iy) y^{s} dy/y
$$
is valid for all $s \in\C$. Analogous to the example in~[\KZ], we derive the
relation
$$
L(g_{18},\chi_D,9) = {2\over\Gamma(9)} (2\pi/|D|)^9 \sum_{n=1}^{\infty} \Kronecker{D}{n} c(D) \phi_8(2\pi n/|D|) \eqno(5.2)
$$
for $g_{18} = \sum_n c(n)n^{-s}$. Here, we write
$$
\phi_8(x) = \int_{1}^\infty y^8 \exp(-xy) dy = {8! \over x^6} \exp(-x) \left(1+x+x^2/2!+\ldots+x^8/8!\right).
$$
The right hand side of~(5.2) converges exponentially fast, and since
we know the Fourier
coefficients of~$g_{18}$ we easily compute the first bound of the Lemma.

Since the space of weight 22 cusp forms is also one-dimensional, the bound for
$K_{12}$ follows analogously. \endproof

\proclaim{Lemma 5.5}{For every $\varepsilon>0$, the twisted $L$-series
associated to the cusp forms $g_{18}$ and $g_{22}$ satisfy
$$
|L(g_{18},\chi_D,9)| \leq B(\varepsilon,9) |D|^{0.5+\varepsilon} \qquad
|L(g_{22},\chi_D,11)| \leq B(\varepsilon,11) |D|^{0.5+\varepsilon}
$$
for all discriminants~$D<0$. Here, $B$ is defined by
$$
B(\varepsilon,n) = {1 \over \sqrt{2\pi}} \max\left\{ \zeta(1+\varepsilon)^2,
\zeta(1+\varepsilon)^2 {\Gamma(n+1/2+\varepsilon) \over \Gamma(n-1/2-\varepsilon)} \right\}.
$$
}
\endproclaim
\noindent
{\bf Proof.} Let $g = \sum_m a(m)q^m$ be either $g_{18}$ or $g_{22}$, and
let $2w$ be the
weight of~$g$. With $\Lambda(s,\chi_D) = \Kronecker{D}{2\pi}^s \Gamma(s+w-1/2)
L(s+w-1/2,g,\chi_D)$, the twisted $L$-series for~$g$ satisfies the functional
equation
$$
\Lambda(s,\chi_D) = \Lambda(1-s,\chi_D)
$$
for all $s \in \C$. We will bound $L(s,g,\chi_D)$ on a vertical line
to the right of the critical strip, which by the functional equation gives
a bound on a vertical line to the left of the critical strip. A variant
of the Phragmen-Lindel\"of theorem will then give the result.

We put $P(s) = \Kronecker{D}{2\pi}^s L(s+w-1/2,g,\chi_D)$ and
$A(m) = c(m)/m^{w-1/2}$. We have $L(s+w-1/2,g_{18},\chi_D) = \sum_m {A(m)
\chi_D(m) \over m^s}$, and the coefficients $A(m)$ are bounded by
$\sigma_0(m) = \sum_{d \mid m} 1$ by Deligne's theorem~[\Del, Th.\ 8.2]. For any
$\varepsilon > 0$ and any $t \in \R$, we bound
$$
|L(1+\varepsilon + w-1/2,g_{18},\chi_D)| \leq \sum_m {|A(m)| \over m^{1+
\varepsilon}} \leq \sum_m {\sigma_0(m) \over m^{1+\varepsilon}} =
\zeta(1+\varepsilon)^2.
$$
We get $|P(1+\varepsilon+it)| \leq {|D|^{1+\varepsilon} \zeta(1+\varepsilon)^2
\over 2\pi} = C_1(\varepsilon) |D|^{1+\varepsilon}$. Using the functional
equation, we bound
$$
|P(-\varepsilon+it)| = \left| {\Gamma(1+\varepsilon-it+w-1/2) \over
\Gamma(-\varepsilon +it+w-1/2)} \zeta(1+\varepsilon -it)^2 \right|
$$
$$
\leq
C_2(\varepsilon) (1+|t|)^{1+2\varepsilon} C_1(\varepsilon) |D|^{1+\varepsilon},
$$
where the last inequality follows from Stirling's formula.

By the Phragmen-Lindel\"of theorem, see e.g.~[\Con, Sec.\ VI.4], we can bound
$$
|P(\sigma+it)| \leq C(\varepsilon) (1+|t|)^{M(\sigma)} |D|^{1+\varepsilon} \qquad\qquad \hbox{for all\ } \sigma\in[-\varepsilon,1+\varepsilon],
$$
where $C(\varepsilon) = \max\{C_1(\varepsilon),
C_2(\varepsilon)C_1(\varepsilon)\}$ is the maximum of the two
$\varepsilon$-dependent bounds on the vertical lines, and $M(\sigma) =
1+\varepsilon-\sigma$ takes the values $M(-\varepsilon) = 1+2\varepsilon$
and $M(1+\varepsilon) = 0$. Taking $\sigma = 1/2$ and $t=0$, we derive
$$
|L(w,g,\chi_D)| \leq \sqrt{2\pi} C(\varepsilon) |D|^{1/2+\varepsilon},
$$
which yields the lemma.
\endproof

\proclaim{Corollary 5.6}{Let the notation be as in Lemma~5.5. Then, for
every $\varepsilon>0$, the coefficients $c_{10,1}(n)$ and $c_{12,1}(n)$ of
$K_{10}$ and $K_{12}$ satisfy
$$
c_{10,1}(n) \leq 320 \sqrt{B(\varepsilon,9)} n^{4.5+1/2\varepsilon}
\qquad \hbox{and} \qquad c_{12,1} \leq 3843 \sqrt{B(\varepsilon,11)}
n^{5.5+1/2\varepsilon}
$$
for all $n \geq 1$.}
\endproclaim
\noindent
{\bf Proof.} Substitute Lemmas 5.4 and 5.5 into Waldspurger's formula~(5.1).
\endproof\par
\ \par
\proclaim{Remark 5.7}{The only room for improvement in Lemma 5.5, and hence
in Corollary 5.6, is in the use of the Phragmen-Lindel\"of theorem. This
theorem yields a factor $n^{1/4}$ in the bound. We can use stronger results
to lower the exponent, but it is harder to make the constants explicit. If
the Lindel\"of hypothesis is true, then the factor $n^{1/4}$ can be
replaced by $n^{\varepsilon}$ (but the constant could in theory be not
explicitly computable).}
\endproclaim

\proclaim{Theorem 5.8}{Let the function~$B$ be as in Lemma~5.5, and
define
$B_2(x) = \exp(2^{1/x}/(x\log 2))$. Then, for every $\varepsilon>0$ and
any $\eta>0$, the
Fourier coefficients $a_{10}(T)$ and
$a_{12}(T)$ of $\chi_{10}$ and $\chi_{12}$ satisfy
\medskip
\halign{\qquad\qquad $#$ & $\leq #$ \cr
|a_{10}(T)| & 320 B_2(\eta) \sqrt{B(\varepsilon,9)} (4\det T)^{4.5+1/2\varepsilon+\eta}\cr
\noalign{\smallskip}
|a_{12}(T)| & 3843 B_2(\eta) \sqrt{B(\varepsilon,11)} (4\det T)^{5.5+ 1/2\varepsilon+\eta}.\cr}
}\endproclaim\noindent
{\bf Proof.} The Fourier coefficient of $\chi_{10}$ for the matrix $T$
is bounded by
$$
\sum_{d \mid (4\det T)} d^9 c_{10,1}\left( {4 \det T \over d^2} \right)
\leq 320\sqrt{B(\varepsilon,9)} \sum_{d \mid (4 \det T)} {(4 \det T)^{4.5+
\varepsilon/2} \over d^{\varepsilon}}.
$$
The sum on the right hand side is bounded by
$$
{(4 \det T)^{4.5+\varepsilon/2}}\sum_{d \mid (4 \det T)} 1 \leq
B_2(\eta) (4\det T)^{4.5+\varepsilon/2+\eta},
$$
see e.g.\ [\HW, Sec.~18.1.]. The proof for $\chi_{12}$ is
similar.
\endproof
\noindent\par
\ \par
\proclaim{Remark 5.9}{If the Lindel\"of-hypothesis is true, then
we get a bound $|a(T)| = O(n^{w/2-3/4+\varepsilon})$ for the Fourier
coefficients of $\chi_{10}$ and $\chi_{12}$. This bound is optimal in the
sense of the Resnikoff-Salda\~na conjecture [\RS]. }
\endproclaim

\head 6. Speed of convergence
\endhead
\noindent
In section we carefully analyse the speed of convergence of the Siegel
Eisenstein series occuring in~(1.2), and this will yield Theorem~1.2 without
too much effort. To analyse the convergence of a Siegel modular function we
a priori have to consider {\it three\/} variables. We begin by showing that
it suffices to look at a `one-dimensional' convergence problem.

The imaginary part $\Im(\tau)$ of a matrix $\tau\in\H_2$ is positive
definite. Hence, there exists $\delta \in \R_{>0}$ with $\Im(\tau) \geq
\delta 1_2$, meaning that $\Im(\tau)-\delta 1_2$ is positive
semi-definite. We define
$$
\delta(\tau) = \sup \{ \delta'\in\R \mid \Im(\tau) \geq \delta' 1_2 \}
$$
to be the `largest' of all these values. With this notation, we have
the following lemma.

\proclaim{Lemma 6.1}{Let $T = \Bigl( {a \atop b/2} \thinspace {b/2 \atop c}
\Bigr) \in \Mat_2({1 \over 4}\Z)$ be positive semi-definite and let
$\tau\in\H_2$. Then the inequality
$$
|\exp(2 \pi i \Tr(T\tau))| \leq \exp(-2 \pi \Tr(T) \delta(\tau))
$$
holds.
}\endproclaim\noindent
{\bf Proof.} We have an equality $|\exp(2 \pi i \Tr(T\tau)) | =
\exp(-2\pi \Tr(T\Im(\tau)))$. Since $T$ is positive semi-definite, we
have $T \Im(\tau) \geq T \delta(\tau)$. The lemma follows.\endproof
\ \par\noindent
We have
$$
E_w(\tau) = \sum_T a(T) \exp(2\pi i \Tr(T\tau)) = \sum_{t=0}^{\infty}
\sum_{T \in S(t)} a(T) \exp(2 \pi i \Tr(T\tau)) \leqno(6.1)
$$
where $S(t)$ is the set of all $2\times 2$ symmetric matrices of
trace~$t$ with
non-negative integer entries on the diagonal and half-integer entries on the
off-diagonal.  The set $S(t)$ clearly has at most $2(t+1)^2$
elements for which $a(T)$ is non-zero.\par
\ \par\noindent
The technique of `splitting up' the evaluation of a Siegel modular form as
in equation~(6.1) enables us to find a lower bound for $|\chi_{10}(\tau)|$.
The idea is that if we have
$$
\Biggl| \sum_{T \in S(t) \atop{\scriptscriptstyle t \leq B}} a(T)
\exp(2\pi i \Tr(T\tau))\Biggr|  > 10 \Biggl| \sum_{T \in S(t) \atop
{\scriptscriptstyle t > B}} a(T) \exp(2\pi i \Tr(T\tau))\Biggr| \eqno(6.2)
$$
then the value of $|\chi_{10}(\tau)|$ is roughly equal to the left hand
side of (6.2). Furthermore, we can apply the upper bound for the Fourier
coefficients of $\chi_{10}$ given by Theorem~5.8 to bound the right
hand side of~(6.2). Taking $B=2$ yields the following lemma.
\ \par\noindent
\proclaim{Lemma 6.2}{Let
$$
M_1 = \Biggl( {1 \atop 0} \thinspace {0 \atop 1} \Biggr), \qquad M_2 = \Biggl(
{1 \atop {1\over 2}} \thinspace {{1\over 2} \atop 1} \Biggr), \qquad M_3 =
\Biggl( {1 \atop -{1\over 2}} \thinspace {-{1\over 2} \atop 1} \Biggr),
$$
and for $\varepsilon,\eta>0$, put
$
M(\varepsilon,\eta) = 320 B_2(\eta) \sqrt{B(\varepsilon,9)},
$
where the notation is as in Theorem~5.8. If, for any $\varepsilon, \eta>0$,
we have $|c| \geq 10 \int_2^\infty 2 M(\varepsilon,\eta) t^{11+\varepsilon+2\eta}
\exp(-2\pi t \delta(\tau)) dt$ for
$$
c = {1\over 2}\exp(2\pi i \Tr(M_1\tau))-{1\over 4}\exp(2\pi i
\Tr(M_2\tau)) -{1\over 4} \exp(2\pi i \Tr(M_3\tau)),
$$
then we have $|\chi_{10}(\tau)| \geq 9/10 |c|$.}
\endproclaim
\noindent
{\bf Proof.} Since $\chi_{10}$ is a cusp form, there are no matrices
$T \in S(0) \cup S(1)$ for which the Fourier coefficient $a(T)$
of $\chi_{10}$ is non-zero. The only matrices $T \in S(2)$ for which
$a(T)$ is nonzero are the matrices $M_1, M_2, M_3$. These matrices have
Fourier coefficients $1/2,-1/4,-1/4$ respectively. Hence, $c$ equals
the left hand side of~(6.2) with $B=2$.

Using Theorem~5.8, we bound the right hand side of~(6.2) from above by
$$
10 \sum_{t=3}^{\infty} 2 t^2 \Bigl| \max_{T \in S(t)}\bigl\{ a(T)  \exp(-2\pi
\Tr(T) \delta(\tau)) \bigr\} \Bigr|
$$
$$
\leq 20 \int_{2}^{\infty}
M(\varepsilon,\eta) t^{11+\varepsilon+2\eta} \exp(-2\pi t \delta(\tau)) dt,
$$
where we used the `AGM-inequality' $4 \det(T) \leq \Tr(T)^2$. The
lemma follows.\endproof
\ \par\noindent
\proclaim{Remark 6.3}{In Lemma~6.2, we can choose any $\varepsilon$ and
$\eta$. The optimal choice depends on the value of~$\delta(\tau)$.}
\endproclaim
\proclaim{Remark 6.4}{If the condition in Lemma~6.2 does not hold for
any $\varepsilon,\delta$, we can
look at the contribution of all matrices of trace $2$ and $3$. If that
majorates the contribution coming from all matrices of trace $4$ and higher,
we have found a lower bound on~$|\chi_{10}(\tau)|$.}
\endproclaim

\ \par\noindent
{\bf Proof of Theorem 1.2.} The Igusa functions are rational expressions in
the Eisenstein series $E_4$, $E_6$ and the cusp forms $\chi_{10}$ and
$\chi_{12}$. The proof consists
of 2 parts: first we analyse the `loss of precision' that occurs when applying
the
formulas~(1.2). Knowing the precision to which to evaluate the four Siegel
modular forms, we then carefully analyse the speed of convergence of these
series.

Using Corollary~3.6, we bound
$$
|a(T)| \leq 19230\, \Tr(T)^{5} \leqno(5.2)
$$
for a Fourier coefficient of $E_{4}$ in case $\det(T)$ is non-zero. For
$\det(T) = 0$ and $\Tr(T) \not = 0$, inequality~(5.2) holds by Theorem~3.4.
We conclude that $|E_4(\tau)|$ is bounded by
$$
1+\int_{0}^{\infty} 2\cdot 19230t^{5}(t+1)^2 \exp(-2\pi t \delta(\tau)) \hbox{d}t
\leq 1+{80\over \delta(\tau)^{8}} + {144 \over \delta(\tau)^7} + {76 \over
\delta(\tau)^6},
$$
and our assumption $\delta(\tau) \geq 1$ implies that we may bound this
by $302$.
For $E_6(\tau)$ we get the bound $|E_6(\tau)| \approx 1+93 / \delta(\tau)^{10}
\leq 94$.
Using Theorem~5.8 with $\eta=1.37$ and $\varepsilon = 0.28$, we derive the
bounds $|\chi_{10}(\tau)| \leq 3487$ and $|\chi_{12}(\tau)| \leq 361893$
for the cusp forms.

Using these four upper bounds, it is straightforward to check that if we
evaluate all four Siegel modular forms up to $k+22$ decimal digits, then
we know the products $\chi_{12}(\tau)^5$, $E_4(\tau) \chi_{12}(\tau)^3$ and
$E_6(\tau)\chi_{12}(\tau)^2$ occuring in formula~(1.2) up to $k$ decimal
digits precision. Furthermore, we know by assumption that $\chi_{10}(\tau)$
does not equal zero. Let $n\in\Z$ be the smallest~$n$ such that
$|\chi_{10}(\tau)| \geq 10^{-n}$ holds. By dividing by $\chi_{10}(\tau)^6$,
we lose $\max\{0,6n\}$ digits precision. Hence, if we evaluate all the
Siegel modular forms occuring in~(1.2) up to $l = k+\max\{22,6n\}$ digits of
precision, we know the Igusa values $j_1(\tau), j_2(\tau), j_3(\tau)$
up to $k$ decimal digits of precision.

We evaluate the Siegel modular functions $E_4,E_6,\chi_{10},\chi_{12}$ using
the sum~(6.1), truncated to only include matrices whose trace is below some
bound~$B$. It
remains to give a value for $B$ such that the function values are accurate
up to $l$ decimal digits. As the speed of convergence of the four series
involved is slowest for $\chi_{12}$, it suffices to look at this function.
Taking $\eta=1.45$ and $\varepsilon=0.1$, we have
$$
\sum_{T \in S(t) \atop {\scriptscriptstyle t \geq B}} a_{12}(T)
\exp(2\pi i \Tr(T \tau)) \leq \int_{B-1}^\infty 524093 t^{15}
\exp(-2\pi t \delta(\tau)) \hbox{d}t
$$
and if the integral is less than $10^{-l}$ then the contribution coming
from the matrices of trace larger than $B$ do not alter the first $l$
decimal digits. The theorem follows. \endproof

\head 7. Examples
\endhead
\noindent
In this section we illustrate the techniques developed in this paper by
evaluating $j_1(\tau)$ for two choices of~$\tau$. 

\subhead 7.1. Example
\endsubhead
\noindent
We detail the evaluation of the Igusa functions $j_1,j_2,j_3$ at
$$
\tau = \Bigl( {2+5i \atop 13+26i} \thinspace {13+26i \atop 83+141i} \Bigr) \in \H_2
$$
to 500 decimal digits of precision. The Igusa functions are rational
expressions in the
Siegel modular forms $E_4,E_6,\chi_{10}$ and $\chi_{12}$, cf.\ Section~1.
The idea is to simply evaluate these series at $\tau$ to high enough
precision and then apply the formulas~(1.2).

We have the rather low bound $\delta(\tau) \geq 0.15$ in this case. However,
for the purpose of evaluating Igusa functions, we may replace $\tau$ by
an $\Sp_4(\Z)$-equivalent matrix~$\tau'$. It is straightforward to check that
the matrix
$$
\tau' = \Bigl( {5i \atop i} \thinspace {i \atop 6i} \Bigr) = 
\left(
\vbox{\vskip.1cm
\halign{\hfill$#$\ \ & \hfill$#$\ \ & \hfill $#$\ \ &\hfill$#$\cr
1& 0&  -2& -13\cr
 -5&   1&  -3& -18\cr
  0&   0&   1&   5\cr
  0&   0&   0&   1\cr
}\vskip-.6cm}
\right)
(\tau)
$$
lies in the fundamental domain for $\Sp_4(\Z) \backslash \H_2$ as 
e.g.~described in~[\Got]. We have $\delta(\tau') \geq 4.3$.

To bound $|\chi_{10}(\tau')|$ from below,
we apply Lemma~6.2. With the notation of this lemma, we compute
$c \approx -1.28 \cdot 10^{-28}$ and the value of the integral is
roughly 
equal to $2 \cdot 10^{-15}$ for $(\eta,\varepsilon) = (1.5,0.1)$. We see 
that Lemma~6.2 does not apply
directly. However, if we compute the contribution $c'$ coming from all
matrices of at most 4, then we get $c' \approx -1.28 \cdot
10^{-28} \approx c$ but we now have
$$
20 \int_{4}^{\infty} 35557 t^{13.2}\exp(-2\pi t\cdot 4.3) \hbox{d}t
\approx 1.1 \cdot 10^{-34}
$$
We conclude that $|\chi_{10}(\tau')|$ is bounded from below by $1.28 \cdot
10^{-28}$.

The lower bound on $|\chi_{10}(\tau')|$ yields that we lose $6\cdot 28 = 168$
decimal digits of precision in the computation of $j_1(\tau')$. However, we
also easily bound
$|\chi_{12}(\tau')|
\leq 4.37 \cdot 10^{-29}$. Hence, we gain $5\cdot 29 =145$ decimal digits
of precision by multiplying by $\chi_{12}(\tau')^5$. The `net loss' of precision
is therefore only $168-145=23$ decimal digits of precision.

Putting everything together, we need to evaluate the Siegel modular forms
$E_4,E_6,\chi_{10},\chi_{12}$ up to $524$ decimal digits precision to know
the values of the Igusa functions up to $500$ decimal digits precision. The
integral
$$
\int_{B-1}^\infty 451485 t^{15.2} \exp(-8.6\pi t) \hbox{d} t
$$
is less than $10^{-524}$ for $B=49$ and we hence have to consider all matrices
of trace up to~$49$.

To compute the Fourier coefficients for all matrices $\Bigl( {a \atop b/2}
\thinspace {b/2 \atop c}\Bigr)$ of {\it trace\/} at
most $49$, we compute the Fourier coefficients of all matrices satisfying
$4ac-b^2 \leq 2401 = 49^2$, with the
convention that we only take the matrices of trace at most $49$ in the
case of determinant~$0$.  To compute all
the coefficients $a(T)$ for $E_4,E_6,\chi_{10}$ and $\chi_{12}$ we compute
the first $2401$ terms of the power series
$$
\theta = 1 + 2 \sum_{n=1}^{\infty} q^{n^2} \hbox{\ \ \quad and\ \ \quad }
\widetilde\theta = \Bigl(1 + 2 \sum_{n=1}^{\infty}
(-1)^n q^{n^2} \Bigr)^4.
$$
Using Proposition~4.3, we compute the first $2401$ coefficients of the modular
forms $H_4$ and $H_6$:
\medskip
\halign{\quad $#$ \thinspace & $#$ & \thinspace $\displaystyle #$\hfill \cr
{-8\over B_4} H_4 &=& 240 + 13440q^3 + 30240q^4 + 138240q^7 + 181440q^8 + 362880q^{11}+O(q^{12})\cr
\noalign{\smallskip}
{-12\over B_{6}} H_6 &=& -504 + 44352q^3 + 166320q^4 + 2128896q^7 + 3792096q^8 + O(q^{11}).\cr}
\medskip\noindent
The coefficients of the forms $H_i$ are the Fourier coefficients of the Jacobi
Eisenstein series~$E_{i}^J$. Using Proposition~5.3 we compute the first $2401$
coefficients of the modular forms $K_{10}$ and $K_{12}$:
\medskip
\halign{\quad $#$ \thinspace & $#$ & \thinspace $\displaystyle #$\hfill \cr
{-1\over 4}K_{10} &=& -1/4q^3+1/2q^4+4q^7-9q^8-99/4q^{11}+O(q^{12})\cr
\noalign{\smallskip}
{1\over 12}K_{12} &=& 1/12q^3+5/6q^4-22/3q^7-11q^8+425/4q^{11}+O(q^{12}).\cr}
\medskip\noindent
The coefficients of $K_{10}$ and $K_{12}$ are the Fourier coefficients of the
Jacobi cusp forms $\varphi_{10,1}$ and $\varphi_{12,1}$.

Since the 4 Siegel modular forms we are interested in lie in the Maa\ss\
Spezialschar, the Fourier coefficient $a(T)$ of one of them only depends
on the determinant of~$T$ and the greatest common divisor of the entries
of~$T$.  We make an array `encoding' these Fourier coefficients as
follows. For every positive integer $N\leq 2401$, we compute its square free part
$N_0$ and write $N = N_0 f^2$. For every divisor $d \mid f$, we compute
and store the
Fourier coefficient belonging to a matrix $T = \Bigl( {a \atop b/2}
\thinspace {b/2 \atop c} \Bigr)$ with $4ac-b^2 = N$ and $\gcd(a,b,c) = d$.
For $E_4$ and $N=16$ we get
$$
[ 997920, 1239840, 0, 1239840 ]
$$
for instance. For $N=0$ we make a list of all positive integers $d\leq X$ and
store the coefficients for the determinant zero matrices with trace~$d$.

The computations so far were independent of the choice of $\tau =
\Bigl( {\tau_1 \atop z} \thinspace {z \atop \tau_2}\Bigr)\in\H_2$. We let
$q_1 = \exp(2\pi i \tau_1)$, $q_2 = \exp(2\pi i \tau_2)$ and
$q_3 = \exp(2\pi i z)$ be the `Fourier variables' of the entries of~$\tau$.
We
compute and store the values $q_1^0=1,q_1,q_1^2,\ldots,q_1^X$ and likewise
for $q_2$. For $\zeta$ we need to compute both the first $X$ powers of $\zeta$
and $\zeta^{-1}$ because the off-diagonal entries of the matrices can be
negative.

The precision needed for this computation is easily computed. Indeed, the
maximum bound for a Fourier coefficient is roughly $10^{21}$ and occurs
for $\chi_{12}$ and a trace~$49$ matrix. As we
need to recognize the values $a(T) \exp(2\pi i \Tr(T\tau))$ up to $524$ decimal
digits precision, we need to compute $q_1$, $q_2$ and $q_3$ with $524+36 =
560$ decimal digits precision.

After making these 4 lists, we now simply loop over $a=0,\ldots,X$, $c=0,
\ldots,X$ and $b=0,\ldots,\lfloor \sqrt{4ac}\rfloor$ and
for the triples $(a,b,c)$ with $b^2-4ac\leq X$ we compute $\gcd(a,b,c)$ and
look up the Fourier coefficient in the stored array.

We implemented this algorithm in the computer algebra package Magma. We did
not attempt to be as efficient as possible in our implementation. On our
64-bit, 2.1 Ghz computer it took roughly 1 second to compute
$j_1(\tau),j_2(\tau), j_3(\tau)$ up to 500 decimal digits precision. We have
$$
j_1(\tau) = 17399743914575167430246482183.29799\ldots
$$
for instance. The computation of the Fourier coefficients of the
Eisenstein series is negligible: the bottleneck is the `loop' over all
matrices $\Bigl( {a \atop b/2} \thinspace {b/2 \atop c} \Bigr)$ satisfying
$0 \leq a \leq X$, $0\leq c \leq X$,
$|b| \leq \lfloor \sqrt{4ac} \rfloor$. \par
\ \par\noindent
%

\subhead 7.2. CM-example
\endsubhead
\noindent
The evaluation of Igusa functions is a main ingredient in the computation
of Igusa class polynomials, which is in turn used to construct 
e.g.\ hyperelliptic curves with cryptographic properties. We illustrate 
our algorithm by recomputing $j_1(\tau)$ for a small CM-point~$\tau$.

Let $K = \Q(\sqrt{-5 + \sqrt{5}})$ be a quartic CM field. The extension~$K/\Q$
is cyclic and $K$ has class number two. Using [\vW, Algorithm 1], see also
[\Weng, Thm.~3.1], we compute that 
$$
\tau' = \Bigl( {2.4060038200i \atop 0.4595058410i} \thinspace {0.4595058410i \atop 
1.9464979789i} \Bigr)
$$
is an approximation to the matrix $\tau$ representing the abelian 
surface~$\C^2/\Phi(\O_K)$, where $\Phi$ is a CM-type for~$K$. We will work
with a 50 digit approximation to~$\tau$.

As shown in~[\vW], the values $j_i(\tau)$ are in fact integers. Hence, we
only need one digit past the decimal place to recognize them and we take
$k=1$ in Theorem~1.2. The matrix $\tau$ already lies in the fundamental 
domain for $\Sp_4(\Z) \backslash \H_2$, and we have $\delta(\tau) \geq 1.66$.
Just as in the previous example, Lemma~6.2 does not apply directly. Using
Remark~6.4, we compute $c \approx -5.3 \cdot 10^{-12}$, where we include all
matrices of trace up to~$6$. The corresponding integral is roughly equal
to~$1.2 \cdot 10^{-16}$ for $\varepsilon=0.1$ and $\eta=1.45$. We conclude
that we may take~$n=12$ in Theorem~1.2. 

Just as in Example~7.1, we bound $|\chi_{12}(\tau)| \leq 3.1 \cdot 10^{-12}$.
We lose at most $1+6 \cdot 12 - 5 \cdot 12 = 13$ digits of precision, and
we need to know the evaluations of the four Siegel modular forms up to
precision~$10^{-14}$. The integral
$$
\int_{B-1}^\infty 524093 t^{15} \exp(-3.32\pi t) \hbox{d}t 
$$
is less than $10^{-14}$ for $B=9$. To get all matrices of trace at most~$9$,
we take all matrices $\Bigl( {a \atop b/2} \thinspace {b/2 \atop c}\Bigr)$ 
satisfying $4ac-b^2 \leq 9^2 = 81$. We compute
$$
j_1(\tau) = 6202728393749.9999\ldots
$$
which is accurate enough to derive $j_1(\tau) = 6202728393750$. 

In this example, it turns out that we only needed to look at the matrices
with $4ac-b^2 \leq 6$. The fact that our bound of 81 was much
higher can be explained as follows. Firstly, our analysis for the precision
loss is for a worst case scenario and we actually do not lose 14 digits of
precision in this example. Secondly, we use the {\it same bound\/} for 
all the
Fourier coefficients of the matrices of a given trace~$t$, whereas these
coefficients actually vary quite a lot.

\head 8. Acknowledgements
\endhead
\noindent
We thank the referee for detailed comments on an earlier
draft of this paper, and Jeff Hoffstein for helpful discussions.


\Refs

\ref\no\Bai
\by{H. Baier}
\paper{Efficient computation of singular moduli with application in cryptography}
\inbook{Fundamentals of computation theory, Springer Lecture Notes in Computer Science}
\vol{2138}
\yr{2001}
\pages{71--82}
\endref

\ref\no\Co
\by{H. Cohen}
\paper{A course in computational algebraic number theory}
\inbook{Springer Graduate Texts in Mathematics, fourth printing}
\vol{138}
\yr{2000}
\endref

\ref\no\CoA
\by{H. Cohen}
\paper{Number Theory, volume II: analytic and modern tools}
\inbook{Springer Graduate Texts in Mathematics}
\vol{240}
\yr{2007}
\endref

\ref\no\CoB
\by{H. Cohen}
\paper{Sums involving the values at negative integers of $L$-functions of quadratic characters}
\jour{Math. Ann.}
\vol{217}
\pages{271--285}
\yr{1975}
\endref

\ref\no\Con
\by{J. B. Conway}
\paper{Functions of one Complex Variable I, 2nd edition}
\inbook{Springer Graudate Texts in Mathematics}
\vol{11}
\yr{1978}
\endref

\ref\no\Del
\by{P. Deligne}
\paper{La conjecture de Weil: I}
\jour{Publ. Math. de l'IHÉS}
\vol{43}
\pages{273--307}
\yr{1974}
\endref

\ref\no\EZ
\by{M. Eichler, D. Zagier}
\paper{The theory of Jacobi forms}
\inbook{Birkh\"auser, Progress in mathematics}
\vol{55}
\yr{1985}
\endref

\ref\no\Got
\by{E. Gottschling}
\paper{Explizite Bestimmung der Randfl\"achen des Fundamentalbereiches der Modulgruppe zweiten Grades}
\jour{Math. Ann.}
\vol{138}
\yr{1959}
\pages{103--124}
\endref

\ref\no\HW
\by{G. H. Hardy, E. M. Wright}
\paper{An introduction to the theory of numbers}
\inbook{Oxford University Press}
\yr{1938}
\endref

\ref\no\IgB
\by{J.-I. Igusa}
\paper{Modular forms and projective invariants}
\jour{Amer. J. Math.}
\vol{89}
\yr{1967}
\pages{817--855}
\endref

\ref\no\IgA
\by{J.-I. Igusa}
\paper{On Siegel modular forms of genus two}
\jour{Amer. J. Math.}
\vol{84}
\yr{1962}
\pages{175--200}
\endref

\ref\no\IgC
\by{J.-I. Igusa}
\paper{On Siegel modular forms of genus two, II}
\jour{Amer. J. Math.}
\vol{86}
\yr{1964}
\pages{392--412}
\endref

\ref\no\Kob
\by{N. Koblitz}
\paper{Introduction to elliptic curves and modular forms}
\jour{Springer Graduate Texts in Mathematics, 2nd edition}
\vol{87}
\yr{1993}
\endref

\ref\no\Koe
\by{M. Koecher}
\paper{Zur Theorie der Modulfunktionen $n$-ten Grades, I}
\jour{Math. Z.}
\vol{59}
\yr{1954}
\pages{399--416}
\endref

\ref\no\Koh
\by{W. Kohnen}
\paper{Estimates for Fourier coefficients of Siegel modular cusp forms of degree 2. II}
\jour{Nagoya Math. J.}
\vol{128}
\yr{1992}
\pages{171--176}
\endref

\ref\no\KZ
\by{W. Kohnen, D. Zagier}
\paper{Values of L-series of Modular Forms at the Center of the Critical Strip}
\jour{Invent. Math}
\vol{64}
\yr{1981}
\pages{175--198}
\endref

\ref\no\Len
\by{H. W. Lenstra Jr.}
\paper{Factoring integers with elliptic curves}
\jour{Ann. of Math. (2)}
\vol{126}
\yr{1987}
\pages{649--673}
\endref

\ref\no\LenPom
\by{H. W. Lenstra Jr., C. Pomerance}
\paper{A rigorous bound for factorig integers}
\jour{J. Amer. Math. Soc.}
\vol{5}
\yr{1992}
\pages{483--516}
\endref

\ref\no\Mah
\by{K. Mahler}
\paper{On a class of non-linear functional equations connected with modular functions}
\inbook{J. Austral. Math. Soc. Ser. A}
\vol{22}
\yr{1976}
\pages{65--118}
\endref

\ref\no\RS
\by{H. L. Resnikoff, R. L. Salda\~na}
\paper{Some properties of Fourier coefficients of Eisenstein series of degree two}
\jour{J. Reine Angew. Math.}
\vol{265}
\pages{90--109}
\yr{1974}
\endref

\ref\no\Ser
\by{J.-P. Serre}
\paper{A course in arithmetic}
\inbook{Springer Graduate Texts in Mathematics}
\vol{7}
\yr{1973}
\endref

\ref\no\Shi
\by{G. Shimura}
\paper{Abelian varieties with complex multiplication and modular functions}
\inbook{Princeton University Press, revised edition}
\yr{1998}
\endref

\ref\no\vW
\by{P. van Wamelen}
\paper{Examples of genus two CM curves defined over the rationals}
\jour{Math. Comp.}
\vol{68}
\yr{1999}
\pages{307--320}
\endref

\ref\no\Wash
\by{L. C. Washington}
\paper{Introduction to cyclotomic fields}
\inbook{Springer Graduate Texts in Mathematics, 2nd edition}
\vol{83}
\yr{1997}
\endref

\ref\no\Weng
\by{A. Weng}
\paper{Constructing hyperelliptic curves of genus 2 suitable for cryptography}
\jour{Math. Comp.}
\vol{72}
\pages{435--458}
\yr{2002}
\endref


\endRefs
\enddocument